\documentclass[12pt,draft]{article}

\usepackage{amsmath}

\begin{document}

\begin{center}
\LARGE\noindent\textbf{A Note on  Hamiltonian Cycles in Digraphs with Large Degrees }\\

\end{center}
\begin{center}
\noindent\textbf{Samvel Kh. Darbinyan }\\

Institute for Informatics and Automation Problems of NAS RA 

E-mail: samdarbin@iiap.sci.am\\
\end{center}

\textbf{Abstract}

 In this note we prove:
 {\it Let $D$ be a 2-strong  digraph of order $n$ such that its $n-1$ vertices have degrees at least $n+k$ and the remaining vertex $z$ has degree at least  $n-k-4$, where $k$ is a non-negative  integer. 
  If $D$ contains a cycle of length  at least $n-k-2$ passing through $z$, then $D$  is Hamiltonian}.

 
\section  {Introduction} 

In this paper, we consider finite digraphs (directed graphs) without loops and multiple arcs. The order of a digraph $D$ is the number of its vertices. We shall assume that the reader is familiar with the standard  terminology on digraphs.  Terminology and  notations not described below follow \cite {[1]}. Every cycle and path is assumed to be simple and directed. 
A cycle (path) in a digraph $D$ is called {\it Hamiltonian} ({\it Hamiltonian path}) if it  includes every  vertex of $D$. A digraph $D$ is {\it Hamiltonian}  if it contains a  Hamiltonian cycle. 

 There are numerous sufficient conditions for the existence of a Hamiltonian cycle in a digraph (see, \cite{[1]},   \cite{[2]}, \cite{[3]}).  The following two   sufficient conditions on the existence of  Hamiltonian cycles in digraphs are classical and famous.

 \textbf{Theorem 1:} (Ghouila-Houri \cite{[4]}). {\it  Let $D$ be a strong digraph of order $n\geq 2$. If for every vertex $x\in V(D)$, $d(x)\geq n$, then $D$ is Hamiltonian}.\\

 \textbf{Theorem 2:} (Meyniel \cite{[5]}). {\it Let $D$ be  a strong digraph of order $n\geq 2$. If $d(x)+d(y)\geq 2n-1$ for all pairs of  non-adjacent  vertices  $x$ and $y$ in $D$, then $D$ is Hamiltonian.}\\

 Nash-Williams \cite{[6]} raised the problem of describing all the extreme  digraphs for Ghouila -Houri's theorem the strong non-Hamiltonian digraphs of order $n$ and with  minimum degree $n-1$. As a solution to this problem, Thomassen \cite{[7]} proved a structural theorem on the extremal digraphs. An analogous problem for the Meyniel theorem (Theorem 2) was considered by the author \cite{[8]}, proving a  structural theorem on the strong non-Hamiltonian digraphs $D$ of order $n$, with the degree condition that $d(x)+d(y)\geq 2n-2$ for every pair of non-adjacent distinct vertices $x, y$. This improves the corresponding structural theorem of Thomassen.  Moreover, in \cite{[8]}, it was also proved that if $m$ is the length of a longest cycle in $D$, then $D$ contains cycles of all lengths $k=2,3,\ldots , m$.

Goldberg, Levitskaya and Satanovskyi \cite{[9]} relaxed the conditions of the Ghouila-Houri theorem by proving  the following theorem.

 \textbf{Theorem 3:} (Goldberg et al. \cite{[9]}). {\it  Let $D$ be a strong digraph of order $n\geq 2$. If $n-1$ vertices of $D$ have degrees at least $n$ and the remaining    vertex has  degree at least $n-1$,  then $D$ is Hamiltonian}.\\
 
 Note that Theorem 3 is  
the best possible in the sense that for every $n$, there is a non-Hamiltonian strong digraphs of order $n$ such that
 its  $n-2$ (or $n-1$) vertices have degrees equal to  $n+1$ (respectively, $n$) and the other two remaining vertices (respectively, the remaining vertex) have degrees equal to $n-1$ (respectively, has degree equal to $n-2$).

It is worth to mention that,  Thomassen \cite{[7]} constructed a strong  non-Hamiltonian digraph of order $n$ with only two vertices of degree $n-1$ and all other vertices have degree not less than $(3n-5)/2$. In \cite{[10]}, it was showed that for every $n\geq 8$ there is a non-Hamiltonian 2-strong digraph of order $n$ such that its $n-1$ vertices have degrees at least $n$ and the remaining vertex has degree 4. 

Taking into account the arguments given above, we can pose the following problem.

\textbf{Problem 1:} Let $D$ be a 2-strong digraph of order $n$ such that its $n-1$ vertices have degrees at least $n$ and the renaining vertex has degree at least $n-k$, where $5\leq n-k\leq n-2$. Investigate the Hamiltonicity of $D$ depending on the values of $n$ and $k$.\\

In \cite{[11]}, it  was reported  that the following theorem holds.

 \textbf{Theorem 4:} (Darbinyan \cite{[11]}). {\it  Let $D$ be a 2-strong  digraph of order $n\geq 9$ with minimum degree at least $n-4$. If $n-1$ vertices of $D$ have degrees at least $n$,   then $D$ is Hamiltonian}.\\ 

 The proof of the last theorem  has never been published. In \cite{[12]}, we  presented the first part of the proof of Theorem 4, by proving the following theorem.

\textbf{Theorem 5:} (Darbinyan \cite{[12]}). {\it  Let $D$ be a 2-strong  digraph of order $n$. Suppose that  $n-1$ vertices of $D$ have degrees at least $n$ and the remaining vertex $z$ has  degree at least $n-4$. If $D$  contains a cycle of length $n-2$ through $z$,   then $D$ is Hamiltonian}.\\

In \cite{[12]}, we also proposed the following conjecture.

\textbf{Conjecture 1:} {\it Let $D$ be a 2-strong  digraph of order $n$ such that its $n-1$ verteices have degrees at least $n+k$ and the remaining vertex  
$z$ has degree at least $n-k-4$, where $k$ is a non-negative integer. Then $D$ is Hamiltonian}.\\ 

Let us note, that Conjecture 1 is an extention  Ghouila-Houri's theorem for 2-strong digraphs and is a generalization of Theorem 4.
The truth of Conjecture 1 in the case $k=0$
follows from Theorem 4. Resently,  we settled Conjecture 1 for any $k\geq 0$.
  Our goal in this note to present the first part of the proof of 
Conjecture 1 for any $k\geq 1$, which we formulate as Theorem 6. 
The second part of the proof (i.e., the complete proof) of Conjecture 1 for any $k\geq 0$ (in particular, the second part of the proof of Theorem 4)  we will present in the forthcoming paper (see arXiv: 2306.16826).\\

\textbf{Theorem 6:} {\it Let $D$ be a 2-strong digraph of order $n\geq 3$ such that its $n-1$ verteices have degrees at least $n+k$ and the remaining vertex  
$z$ has degree at least $n-k-4$, where $k\geq 0$ is an integer. If $D$ has a cycle of length at least $n-k-2$ through $z$, then $D$ is Hamiltonian}.\\

\section{ Further Terminology and Notation}

 For the sake of clarity we repeat the most impotent definition. The vertex set and the arc set of a digraph $D$ are    denoted 
  by $V(D)$  and   $A(D)$, respectively.   The arc of a digraph $D$ directed from $x$ to $y$ is denoted by $xy$ or $x\rightarrow y$ (we also say that $x$ {\it dominates} $y$ or $y$ is an
  {\it out-neighbour} of $x$ and $x$ is an {\it in-neighbour} of $y$), and $x\leftrightarrow y$ denotes that 
$x\rightarrow y$ and $y\rightarrow x$ ($x\leftrightarrow y$ is called  {\it 2-cycle}).
   If $x\rightarrow y$ and $y\rightarrow z$, we write 
$x\rightarrow y\rightarrow z$.

 Let $A$ and $B$ be two disjoint subsets in $V(D)$. The notation $A\rightarrow B$ means that  every
   vertex of $A$  dominates every vertex of $B$. 
 We define 
$A(A\rightarrow B)=\{xy\in A(D)\, |\, x\in A, y\in B\}$. 
If $x\in V(D)$
   and $A=\{x\}$ we sometimes write $x$ instead of $\{x\}$.
Let $N_D^+(x)$, $N_D^-(x)$ denote the set of  out-neighbors, respectively the set  of in-neighbors of a vertex $x$ in a digraph $D$.  If $A\subseteq V(D)$, then $N_D^+(x,A)= A \cap N_D^+(x)$ and $N_D^-(x,A)=A\cap N_D^-(x)$. 
The {\it out-degree} of $x$ is $d_D^+(x)=|N_D^+(x)|$ and $d_D^-(x)=|N_D^-(x)|$ is the {\it in-degree} of $x$. Similarly, $d_D^+(x,A)=|N_D^+(x,A)|$ and $d_D^-(x,A)=|N_D^-(x,A)|$. The {\it degree} of the vertex $x$ in $D$ is defined as $d_D(x)=d_D^+(x)+d_D^-(x)$ (similarly, $d_D(x,A)=d_D^+(x,A)+d_D^-(x,A)$). We omit the subscript if the digraph is clear from the context. The subdigraph of $D$ induced by a subset $A$ of $V(D)$ is denoted by $D\langle A\rangle$. In particular, $D-A=D\langle V(D)\setminus  A\rangle$.
For integers $a$ and $b$, $a\leq b$, by $[a,b]$ we denote the set $\{x_a,x_{a+1},\ldots , x_b\}$. If $j<i$, then $\{x_i,x_{i+1},\ldots , x_j\}=\emptyset$. 
 A path is a digraph with vertex set   $\{x_1,x_2,\ldots ,x_k\}$  and  arc set  $\{x_1x_{2}, x_2x_3,\ldots , x_{k-1}x_{k}\}$, and is denoted by  $x_1x_2\cdots x_k$. This is also called an $(x_1,x_k)$-path or a path from $x_1$ and $x_k$. If we add the arc $x_kx_1$ to the path above, we obtain a cycle  $x_1x_2\ldots x_kx_1$. The {\it length} of a cycle or a path is the number of its arcs.  If $P$ is a path containing a subpath from $x$ to $y$, we let $P[x,y]$ denote that subpath. Similarly, if $C$ is a cycle containing vertices $x$ and $y$, $C[x,y]$ denotes the subpath of $C$ from $x$ to $y$, and an $(x,y)$-path $P$ is a $C$-bypass (or is a $(C,x,y)$-bypass) if $|V(P)|\geq 3$ and $V(P)\cap V(C)=\{x,y\}$.   
Let $D$ be a digraph and $z\in V(D)$. By $C_m(z)$ (respectively, $C(z)$) we denote  a cycle in $D$ of length $m$ (respectively, any cycle in $D$), which contains the vertex   $z$. Similarly, we denote by $C_k$ a cycle of length $k$.
  A digraph
 $D$ is {\it strong} ({\it strongly connected}) if, for every pair $x, y$ of distinct vertices in $D$, there exists an $(x,y)$-path and a $(y,x)$-path.  A digraph
 $D$ is {\it $k$-strong} ({\it $k$-strongly connected}) if, $|V(D)|\geq k+1$ and for any set $A$ of at most $k-1$ vertices $D-A$ is strong.
Two distinct vertices $x$ and $y$ are {\it adjacent} if $xy\in $ or $yx\in V(D) $ (or both). The {\it converse}  digraph of  $D$ is the digraph obtained from  $D$ by reversing the direction of all  arcs in $D$. 
 We will use {\it 
 the principle of digraph duality}: Let $D$ be a digraph, then $D$ has a subdigraph $H$ if and only if the converse digraph of $D$ has the converse of the subdigraph $H$.

\section {Preliminaries}
    
 In our proofs, we will use the following well-known simple lemmas.
    
 \textbf{Lemma 1:} (H\"aggkvist and Thomassen \cite{[13]}). {\it Let $D$ be a digraph of order $n\geq 3$ containing a cycle $C_m$ of length $m$, $m\in [2,n-1]$. Let $x$ be a vertex not contained in this cycle. If $d(x,V(C_m))\geq m+1$, then for every $k\in [2,m+1]$, $D$ contains a cycle  $C_k$ of length $k$ including $x$.} 
  
  The next lemma is a slight  modification of the lemma by Bondy and Thomassen \cite{[14]}, it is very useful and will be used extensively throughout this paper. 
  
 \textbf{Lemma 2:}. {\it Let $D$ be a digraph of order $n\geq 3$ containing a path $P:=x_1x_2\ldots x_m$, $m\in [2,n-1]$. Let $x$ be a vertex not contained in this path. If one of the following condition holds: 
 
 (i) $d(x,V(P))\geq m+2$, 
 
 (ii) $d(x,V(P))\geq m+1$ and $xx_1\notin A(D)$ or $x_mx\notin A(D)$, 
 
 (iii) $d(x,V(P))\geq m$ and $xx_1\notin A(D)$ and $x_mx\notin A(D)$,
 
   then there is an $i\in [1,m-1]$ such that $x_i\rightarrow x\rightarrow x_{i+1}$, i.e.,  $D$ contains a path $x_1x_2\ldots x_ixx_{i+1}\ldots x_m$ of length $m$ (we say that $x$ can be inserted into $P$)}.  \\

 In \cite{[10]}, the author proved the following theorem.
 
 \textbf{Theorem 7:} (Darbinyan  \cite{[10]}). {\it Let $D$ be a strong digraph of order $n\geq 3$. Suppouse that $d(x)+d(y)\geq 2n-1$ for all pairs of  non-adjacent  vertices  $x, y\in V(D)\setminus \{z\}$, where $z$ is an arbitrary vertex in $V(D)$. Then $D$ is Hamiltonian or contains a cycle of length $n-1$.}\\
  
  Using Theorem 7 and Lemma 1,  it is not difficult to show that the following corollary  is true.

  \textbf{Corollary 1:} {\it Let $D$ be a strong digraph of order $n\geq 3$. Suppose that $n-1$ vertices of $D$ have degrees at least $n$.  Then $D$ is Hamiltonian or contains a cycle of length $n-1$ (in fact, $D$ has a cycle that contains all the vertices of degrees at least $n$)}.\\ 
   
By the same arguments as the proof Lemma 4 in \cite{[12]}, we can prove the following lemma.

\textbf{Lemma 3:} {\it Let $D$ be a non-Hamiltonian digraph of order $n\geq 4$ such that its  $n-1$ vertices  have degrees at least $n$ and the remaining vertex $z$ has degree at most  $n-2$.
Suppose that $C_m(z)=x_1x_2\ldots x_mx_1$ is a longest cycle of length $m$ through $z$ in $D$. If  $D$    has a $(C_m(z), x_i,x_j)$-bypass, then  $z\in V(C_m(z)[x_{i+1},x_{j-1}])$}.

 \section*{Proof of Theorem 6}

  Our proofs are based on the arguments of \cite{[12]}.
 First, we will prove the following lemma, which is of independemt interest.

 \textbf{Lemma 4:} {\it Let $D$ be a  non-Hamiltonian 2-strong  digraph of order $n$  such that its $n-1$ vertices have degrees at least $n$ and the remaining    vertex $z$ has degree at most $n-2$.
 Suppose that $C_{m+1}(z)=x_1x_2\ldots x_mzx_1$ with  $m\in [2,n-3]$ is a longest cycle  through $z$ in $D$. If two distinct 
 vertices $y_1$, $y_2$  of $Y:= V(D)\setminus V(C_{m+1}(z))$ are  mutually reachable in $D\langle Y\rangle$ and   for each $y_i\in \{y_1,y_2\}$,
 $d(y_i,\{x_1,x_2,\ldots , x_m\})=m+1$, then $n\geq 6$ and  $d(z,\{x_1,x_2,\ldots , x_m\})\leq m-2$.}

 \textbf{Proof:} By contradiction, suppose that $d(z,\{x_1,x_2,\ldots , x_m\})\geq m-1$. By $P$ we denote the path $x_1x_2\ldots x_m$. It is clear that   $|Y|=n-m-1$.  Since  $C_{m+1}(z)$ is a longest cycle, it follows that every vertex $y\in Y$ cannot be inserted into $C_{m+1}(z)$. Then  by Lemma 1, 
$d(y,V(C_{m+1}(z)))\leq m+1$ and
$$
n\leq d(y)=d(y,V(C_{m+1}(z)))+d(y,Y)\leq m+1+d(y,Y).    
$$
Hence, $d(y,Y)\geq n-m-1=|Y|$. Therefore by Ghoula-Houri's theorem, $D\langle Y\rangle$ contains a Hamiltonian path. Let $H_1$, $H_2, \cdots , H_f$ be the strong components of
$D\langle Y\rangle$ labelled in such way that no vertex of  $H_i$ dominates a vertex of $H_j$  whenever $i>j$. Since $D\langle Y\rangle$  has a Hamiltonian path, it follows that for each $i\in [1,f-1]$ there is an arc from $H_i$ to $H_{i+1}$. From $d(y,Y)\geq |Y|$ we obtain that if $u\in V(H_i)$, then $d(u,V(H_{i}))\geq |V(H_i)|$.   Therefore,  $|V(H_i)|\geq 2$ and by Ghoula-Houri's theorem every subdigraph $H_i$ is Hamiltonian. Observe that for some $d\in [1,f]$ the vertices $y_1$, $y_2$ are in $V(H_d)$.
Let $R$ be a longest path between $y_1$ and $y_2$ in $D\langle Y\rangle$.
Since $D$ is 2-strong and $C_{m+1}(z)$ is a longest cycle, using Lemma 2 and the fact that  $d(y_i,V(P))=m+1$ it is not difficult to show   that there is an integer $l\in [2,m-1]$ such that
$$
\{x_{l}, x_{l+1}, \ldots , x_{m}\}\rightarrow \{y_1,y_2\}\rightarrow \{x_{1}, x_2, \ldots , x_l\}.    \eqno(1)
$$
Put
   $E:=\{x_{1}, x_2, \ldots , x_{l-1}\}$  and  $F:=\{x_{l+1}, x_{l+2}, \ldots , x_{m}\}$. Note that $|E|\geq 1$ and $|F|\geq 1$. 
Since $C_{m+1}(z)$ is a longest cycle,  from (1) it follows that 
$$
A(\{z\}\cup E\rightarrow \cup_{i=1}^dV(H_i))=A(\cup_{i=d}^fV(H_i)\rightarrow \{z\}\cup F)=\emptyset,  \eqno (2)
$$
in particular, $d(z,V(H_d))=0$. 
Note that from $|Y|\geq 2$, $|E|\geq 1$ and $|F|\geq 1$   it follows that $n\geq 6$.  

We claim that $d(z,Y)=0$. Assume that this is not the case. Then by (2), for some $u\in \cup_{i=1}^{d-1}V(H_i)$
or $w\in \cup_{i=d+1}^{f}V(H_i)$, $u\rightarrow z$ or $z\rightarrow w$. Without loss of generality assume that $u\rightarrow z$. Then  $A(\{x_1, x_2,\ldots , x_m\}\rightarrow V(H_1))=\emptyset$ for otherwise for some $i\in [1,m]$ there is a $(C_{m+1}(z),x_i,z)$-byypass, a contradiction to Lemma 3. Therefore,  $A(\cup_{i=2}^{f}V(H_i)\cup \{z,x_1, x_2,\ldots , x_m\}\rightarrow V(H_1))=\emptyset$, which contradicts that $D$ is 2-strong.

If $u,v\in Y$, then we use $L(u,v)$ to denote a $(u,v)$-path in $D\langle Y\rangle$.    We need to prove  Claims 1-2 bellow.\\

 \textbf{Claim 1.} 
 
 {\it (i) If $d^-(z,E)\geq 1$, then $d^+(z,F)=0$.  
 (ii)  $A(E\rightarrow F)\not=\emptyset$}.

 \textbf{Proof.} (i) By contradiction, suppose that  $x_i\in E$,  $x_j\in F$ and $x_i\rightarrow z\rightarrow x_j$. Then by (1), $y_1\rightarrow x_{i+1}$ and $x_{j-1}\rightarrow y_2$. Hence, $C_{m+3}(z)= x_1x_2\ldots x_izx_j\ldots x_{m}y_1x_{i+1}\ldots x_{j-1}y_2x_1$, a contradiction. 
   
 (ii)  Again  by contradiction, suppose that 
 $A(E\rightarrow F)=\emptyset$. Let us consider the subdigraph $D-x_l$. Then by (2),  we have: if $d^-(z,E)=0$, then 
$$
A(\cup_{i=d}^fV(H_i)\cup E\rightarrow \cup_{i=1}^{d-1}V(H_i)\cup \{z\}\cup F)=\emptyset,  
$$
and if $d^-(z,E)\geq 1$, then $d^+(z,F)=0$ (Claim 1(i)) and 
$$
A( \cup_{i=d}^fV(H_i)\cup \{z\}\cup E\rightarrow \cup_{i=1}^{d-1}V(H_i)\cup F
)=\emptyset. 
$$
Thus,  in both cases,
  $D-x_l$ is not strong,  which contradicts  that $D$ is 2-strong. \fbox\\

   From now on, we assume that  $x_ax_b\in A (E\rightarrow F)$. Note that by 
(1), we have 
$$
x_{b-1}\rightarrow \{y_1,y_2\}\rightarrow x_{a+1}. \eqno (3)
$$ 
From Claim 1(i) it follows that either $x_a z\notin A(D)$ or $z x_b\notin A(D)$. Therefore, 
   since $z$ cannot be  inserted into $P$, using Lemma 2(ii), we obtain
   $$
   d(z, \{x_1,x_2,\ldots, x_a\})+d(z, \{x_b,x_{b+1},\ldots, x_m\})\leq a+m-b+2.   \eqno (4)
   $$

  It is clear that $|V(R)|\geq 2$ since $H_d$ is Hamiltonian and $|V(H_d)|\geq 2$. Using (3), $|V(R)|\geq 2$ and  the fact that $C_{m+1}(z)$ is a longest cycle, it is not difficult to prove that the following claim is true.

  \textbf{Claim 2.} 
  
  (i) {\it If $i\in [a+1,l-1]$, then $x_{i}z\notin A(D)$.
   
    (ii) If $j\in [l+1,b-1]$, then $zx_{j}\notin A(D)$.
    
    (iii) If $i\in [a+1,l]$  and $i-a\leq |V(R)|$, then  $zx_{i}\notin A(D)$.
    
    (iv) If $j\in [l,b-1]$ and $b-j\leq |V(R)|$, then $x_{j}z\notin A(D)$}. 

  \textbf{Proof.} By contradiction: (i) Assume that $i\in [a+1,l-1]$ and $x_i\rightarrow z$. Then by (2) and (3), we have $y_1\rightarrow x_{i+1}$, $x_{b-1}\rightarrow y_2$ and $y_2\rightarrow x_{a+1}$. Therefore, $C_{m+3}(z) =x_1x_2\ldots x_ax_b\ldots x_my_1\\x_{i+1}\ldots x_{b-1}y_2x_{a+1}\ldots x_izx_1$, a contradiction.

(iii) Assume that $i\in [a+1,l]$, $i-a\leq |V(R)|$ and $z\rightarrow x_i$. Then the cycle $C(z)=x_1x_2\ldots x_ax_b\ldots x_mzx_{i}\ldots x_{b-1}Rx_{1}$ has length at least $m+2$, a contradiction.

Similarly, we can prove that (ii) and (iv) also are true. \fbox

 Now we will consider the following cases depending on the values of $a$ and $b$ with respect to $l$.

  \textbf{Case 1.} $a\leq l-3$ and $b\geq l+2$.

  Then by Claim 2,
   $d(z,\{x_{a+1},x_{a+2},x_{b-1}\})=0$ and $x_{b-2}z\notin A(D)$. Therefore, since $z$ cannot be inserted into $P$, using (4) and Lemma 2, we obtain
  $$
  m-1\leq  d(z,V(P))=d(z, \{x_1,x_2,\ldots, x_a,x_b,x_{b+1},\ldots, x_m\})$$ $$+d(z,\{x_{a+3},x_{a+4},\ldots , x_{b-2}\})\leq a+m-b+2+b-2-a-2=m-2,  
   $$
  which is a contradiction.

Befor consider Cases 2-4, we will show the following two propositions.

 \textbf{Proposition 1.} If $b=l+1\leq m-1$, then 
$$A(\{x_1,x_2,\ldots , x_l\}\rightarrow \{x_{l+2},x_{l+3},\ldots , x_m\})\not= \emptyset.$$

\textbf{Proof.} Assume that the proposition is  not true. Note that $x_lz\notin A(D)$ by Claim 2.
Since $D-x_{l+1}$ is strong, it follows that in $D-x_{l+1}$ there is a path from a vertex $x_i\in \{x_1,x_2,\ldots , x_l\}$ to a vertex $x_j\in \{x_{l+2},x_{l+3},\ldots , x_m\}$. Let $Q$ be such a shortest  path. Then from (2), $d(z,Y)=0$ and Claim 1(i) it follows that the internal vertices of $Q$ are in $\cup_{i=1}^{d-1}V(H_i)$. 
 This means that $x_i=x_l$ and  $C_{m+1}(z)$ contains an $(x_l,x_j)$-bypass such that $z\notin V(C_{m+1}(z)[x_{l+1},x_{j-1}])$.  Therefore by Lemma 3, $D$ contains a $C(z)$-cycle of length at least $m+2$, which contradicts  that   $C_{m+1}(z)$ is a longest cycle in $D$. \fbox \\\\

\textbf{Proposition 2.} Suppose that $u\in V(H_f)$ and the vertices $u$ and $x_s$ with $s\geq 1$ are adjacent. If  $d(u, \{x_1,x_2,\ldots , x_{s-1}\})=0$, then $|V(H_f)|\geq s+1$. 
 
\textbf{Proof.} Recall that every vertex  of $V(H_f)$  cannot be inserted into $P=x_1x_2\ldots x_m$. Then by Lemma 2, $d(u,V(P))=d(u,\{x_s,x_{s+1},\ldots , x_m\})\leq m-s+2$. Therefore, 
$$
n\leq d(u)= d(u,V(H_f))+d^-(u, \cup_{i=1}^{f-1}V(H_i))+d(u,V(P))$$ 
$$
\leq d(u,V(H_f))+n-m-1-|V(H_f)|+ m-s+2,
$$
i.e., $d(u,V(H_f))\geq |V(H_f)|+s-1$. This together with $d(u,V(H_f))\leq 2|V(H_f)|-2$ implies that $|V(H_f)|\geq s+1$. \fbox \\\\

  \textbf{Case 2.} $a\leq l-3$ and $b= l+1$.
  
Then 
   by Claim 2,
   $d(z,\{x_{a+1},x_{a+2}\})=0$, $x_{b-1}z\notin A(D)$. 
Again using (4) and Lemma 2, we obtain
  $$
  m-1\leq d(z, \{x_1,x_2,\ldots, x_a,x_b,x_{b+1},\ldots, x_m\})+d(z,\{x_{a+3}, x_{a+4},\ldots , x_{b-1}\})$$ $$\leq a+m-b+2+b-1-a-2= m-1.  
   $$
 Therefore, 
$
d(z,\{x_{a+3},x_{a+4},\ldots ,x_{b-1}\})= b-a-3,
$
and hence  by Claim 2, 
$$
z\rightarrow \{x_{a+3},\\x_{a+4}, \ldots , x_l\}. \eqno (5)
$$
Now it is easy to see that $|V(H_d)|=|V(R)|=2$, i.e., $R=y_1y_2$, for otherwise $|R|\geq 3$ and $C(z)=x_1x_2\ldots x_ax_b\ldots x_mzx_{a+3}\ldots x_{b-1}Rx_1$ has length greater than $m+1$, a contradiction.
Notice that from $|V(H_d)|=2$, $d(y_j,V(P))=m+1$, $d(y_j,\{z\})=0$ and $d(y_j)\geq n$ it follows that 
$$
\cup^{d-1}_{i=1}V(H_i)\rightarrow \{y_1,y_2\}\rightarrow \cup^{f}_{i=d+1}V(H_i). \eqno (6)        
$$

\textbf{Subcase 2.1.} $a\leq l-3$, $b=l+1$ and $m\geq l+2$.

Taking into account  Case 1, we may assume that
$$
A(\{x_1,x_2,\ldots , x_{a}\}\rightarrow 
\{x_{l+2},x_{l+3},\ldots , x_{m}\})=\emptyset. 
$$
This together with Proposition 1 implies that there are
$i\in [a+1,l]$ and $j\in [l+2,m]$ such that $x_i \rightarrow x_j$.  Now using (3) and (5), we obtain: If $i\in [a+2,l-1]$, then 
$C_{m+3}(z)=x_1x_2\ldots x_ix_j\ldots x_mzx_{i+1}\ldots x_{j-1}Rx_1$, 
if $i=a+1$, then $C_{m+2}(z)=x_1x_2\ldots x_ix_j\ldots x_mzx_{a+3}\\\ldots x_{j-1}Rx_1$, and if $i=l$, then
$C_{m+3}(z)=x_1x_2\ldots x_ax_{l+1}\ldots  x_{j-1}Rx_{a+1}\ldots x_lx_j\ldots x_mzx_1$. Thus, in either case, we have a $C(z)$-cycle of length at least $m+2$, a contradiction.\\

\textbf{Subcase 2.2.} $a\leq l-3$,  $b=l+1=m$ and $d\geq 2$.

 If $x_{m-1}\rightarrow w$ for some $w\in V(H_1)$, then using (5) and (6), we obtain  that $C_{m+2}=x_1x_2\ldots x_ax_mzx_{a+3}\ldots x_{m-1}wy_1y_2x_1$, a contradiction. We may assume that $d^+(x_{m-1},\\ V(H_1))=0$. This together with (2) and  
$A( \cup^{f}_{i=2}V(H_i)\rightarrow V(H_1))=\emptyset$   implies that  $A(\cup^{f}_{i=2}V(H_i)\cup \{x_1,x_2, \ldots ,x_{m-1},z\}\rightarrow V(H_1))=\emptyset$, which means that $D-x_m$ is not strong, a contradiction.\\

\textbf{Subcase 2.3.} $a\leq l-3$, $b=l+1=m$ and $d=1$, i.e, $V(H_1)=\{y_1,y_2\}$.

  Let $x_a\rightarrow x_m$ and $a$ is the minimum with this property.

If $x_i\rightarrow x_m$ with $i\in [a+1,m-2]$,  then, since $z\rightarrow x_{i+1}$ or  $z\rightarrow x_{i+2}$,  the cycle $C(z)=x_1\ldots x_ix_mzx_{i+1}(or x_{i+2})\ldots x_{m-1}y_1y_2x_1$ has length at least $m+2$, a contradiction. We may therefore assume that 
$$
d^-(x_m,\{x_{a+1},x_{a+2},\ldots , x_{m-2}\})=0. \eqno(7)
$$ 

Let $u\in \cup_{i=2}^fV(H_i)$. If $u\rightarrow x_1$, then by (5) and (6),  $C_{m+2}=x_1x_2\ldots x_ax_mzx_{a+3}\ldots \\x_{m-1}y_1y_2ux_1$, a contradictin.  If  $x_1\rightarrow u$, then by Lemma 3, $A(V(H_f)\rightarrow \{x_2,x_3,\ldots ,\\x_m,z \})=\emptyset$. As a result, we have
 $A(V(H_f)\rightarrow\cup_{i=1}^{f-1}V(H_i)\cup \{x_2,x_3,\ldots ,x_m,z \})=\emptyset$, which contradicts that $D-x_1$ is strong. Thus, we may assume that
$$d(x_1, \cup_{i=2}^fV(H_i))=0.  \eqno (8) $$

Assume first that $x_1x_m\notin A(D)$.   Then $a\geq 2$. Since $d^-(z)\geq 2$, $d(z,Y)=0$ and $d^-(z,\{x_{a+1},x_{a+2},\ldots ,x_{m-1}\})=0$, it follows that $d^-(z,\{x_1,x_2,\ldots , x_a\})\geq 1$. This together with Claim 1(i) implies that $zx_m\notin A(D)$. Hence, $d(z,\{x_{a+1},x_{a+2},\ldots , x_m\})\leq m-a-2$. From this and $d(z,V(P))\geq m-1$ we have that $d(z,\{x_1,x_2,\ldots  ,x_a\})\geq a+1$. Therefore, since $z$ cannot be inserted into the path $x_1x_2\ldots x_a$, using Lemma 2, we obtain $x_a\rightarrow z$.

 Let $j\in [a+1,a+2]$.  Then we have: if  $x_1\rightarrow x_j$, then $C(z)=x_1x_j\ldots x_my_1y_2x_2\ldots x_azx_1$ has length at least $m+2$, if  $x_j\rightarrow x_1$, then  $C_{m+2}(z)=x_1x_2\ldots x_azx_{a+3}\ldots x_my_1y_2x_jx_1$. Thus, in either case we have a contradiction, which proves that
$$
d(x_1,\{x_{a+1},x_{a+2}\})=0.     \eqno (9)
$$
From (9) it follows that $f\geq 2$, for otherwise the cycle $C_{n-2}(z)= x_1x_2\ldots x_azx_{a+3}\ldots x_m\\y_1y_2x_1$  does not contain the vertices $x_{a+1}$, $x_{a+2}$ and 
$d(x_1,\{x_{a+1},x_{a+2}\})= d(z,\{x_{a+1},x_{a+2}\})\\=0$, which is impossible.

Now we want to show that 
$$
A(V(H_f)\rightarrow 
\{x_{a+1},x_{a+2},\ldots , x_{m-1}\})=\emptyset. \eqno (10)
$$
Assume that (10) is not true, i.e., there are two vertices $u\in V(H_f)$ and $x_j$ with $j\in [a+1,m-1]$ such that
$u\rightarrow x_j$. Recall that $C(z)=x_1x_2\ldots x_azx_{a+3}\ldots x_my_1y_2x_1$ is a cycle of length $m+1$ and the vertices $x_{a+1}$, $ x_{a+2}$ are not on $C(z)$. If $j\in [a+2,m-1]$, then by Lemma 3, $d^+(x_{a+1}, Y)=0$, and $x_{a+1}x_2\notin A(D)$ since
$y_2\rightarrow x_{a+1}$. Since $x_{a+1}$ cannot be inserted into $C(z)$, using Lemma 2 and the fact that $d(x_{a+1},\{x_1,z\})=0$, we obtain
$$
n\leq d(x_{a+1})= d(x_{a+1},\{x_{a+2}\})+d(x_{a+1}, V(C(z)[x_{a+3},y_2])) $$ $$+d(x_{a+1}, \cup_{i=2}^f(V(H_i))) +d(x_{a+1},\{x_2,x_3,\ldots ,x_a\})$$ $$\leq2+m+1-a-1+1+n-m-1-2+a-1=n-1,
$$
a contradiction. Thus, 
$$
A(V(H_f)\rightarrow 
\{x_{a+2},x_{a+3},\ldots , x_{m}\})=\emptyset. \eqno (11)
$$
If $j=a+1$, i.e., $u\rightarrow x_{a+1}$, then by Lemma 3,
$$
A(\{x_1,x_2,\ldots x_a\}\rightarrow V(H_f))=\emptyset. \eqno(12)
$$
If $A(V(H_f)\rightarrow 
\{x_{1},x_{2},\ldots , x_{a}\})=\emptyset$, then by (11) and (12), we have $$A(V(H_f)\rightarrow \cup_{i=1}^{f-1}V(H_i)\cup 
\{z, x_{1},x_{2},\ldots , x_{a}, x_{a+2},x_{a+3},\ldots ,x_m\})=\emptyset,
$$
which means that $D-x_{a+1}$ is not strong. We may therefore assume that $A(V(H_f)\rightarrow 
\{x_{1},x_{2},\ldots , x_{a}\})\not=\emptyset$. Let $u\rightarrow x_s$, where $u\in V(H_f)$, $s\in [1,a]$ and $s$ is the minimum with this property, i.e.,  $A(V(H_f)\rightarrow 
\{x_{1},x_{2},\ldots , x_{s-1}\})=\emptyset$. Since $d(x_1,V(H_f))=0$ (by (8)), it follows that $s\geq 2$. Then by (12), $A(V(H_f),
\{x_{1},x_{2},\ldots , x_{s-1}\})=\emptyset$.    By Proposition 2, $|V(H_f)|\geq s+1$. Therefore, the cycle $C(z)=y_1y_2H(v,u)x_sx_{s+1}\ldots x_azx_{a+3}\ldots x_my_1$, where $H(v,u)$ is a Hamiltonian path in $H_f$, has length  $m+3$, a contradiction.

Assume second that $x_1\rightarrow x_m$. Then from Claim 2 and $d(z,Y)=0$ it follows that $d^-(z,Y\cup \{x_2,x_3,\ldots ,x_{m-1}\})=0$, which in turn implies that $x_1\rightarrow z$ since $d^-(z)\geq 2$. Then by Claim 1(i), $zx_m\notin A(D)$. By (7), we have that 
$d^-(x_m,\{x_2,x_3,\ldots , x_{m-2},z\}\cup Y)=0$. Since $D-x_{m-1}$ is strong, it contains a path from a vertex $x_j\in \{x_2,x_3,\ldots ,x_{m-2}\}$ to the vertex $x_m$. Let $Q$ be such a shortest  path. Using (7)  and (8), it is not difficult to see that  $Q=x_jx_1x_m$.  Therefore the cycle $C(z)=x_jx_1zx_{j+1}(or x_{j+2}) \ldots x_my_1y_2x_2\ldots x_j$ has length at least $m+2$, which contradicts that a longest cycle through $z$ in $D$ has length $m+1$.

  \textbf{Case 3.} $a=l-2$.  

Taking into account the case $a\leq a-3$ and $b\geq l+2$, we may assume that $b\leq l+2$. 

\textbf{Subcase 3.1.} $b=l+2$. 

Then by Claim 2,
$d(z,\{x_{a+1},x_{a+2}, x_{a+3}=x_{b-1}\})=0$. Therefore by (4), we have 
$$
m-1\leq d(z,V(P))=d(z,\{x_1,x_2,\ldots ,x_a, x_b,x_{b+1}, \ldots ,x_m\})$$ $$\leq a+m-b+2=b-4+m-b-2=m-2,
$$
a contradiction. 

\textbf{Subcase 3.2.} $b=l+1$.

Then $d(z,\{x_{a+1},x_{a+2} \})=0$. Now using (4), we obtain that $d(z,V(P))=m-1$.

Assume first that $m \geq l+2$. 
 Taking into account the considered cases, we may assume that
$$
 A(\{x_1,x_2,\ldots , x_{a}\}\rightarrow \{x_{b+1},x_{b+2},\ldots , x_{m}\})= \emptyset. \eqno (13)
 $$ 

 If $x_i\rightarrow x_j$ with $i\in [a+1,a+2=l]$ and $j\in [l+2,m]$, then the cycle $C(z)=x_1x_2\ldots x_ax_{l+1}\ldots x_{j-1}Rx_ix_j\ldots x_mzx_1$ is a cycle of length at least $m+2$, a contradiction. Therefore, we may assume that $A(\{x_{l-1},x_l\}\rightarrow \{x_{l+2},x_{l+3},\ldots , x_{m}\})=\emptyset$.
   This together with (13) implies that 
$
 A(\{x_1,x_2,\ldots , x_{l}\}\rightarrow 
\{x_{l+2},x_{l+3},\ldots , x_{m}\})=\emptyset 
 $, 
which a contradics Proposition 1.

Assume second that $b=l+1=m$.
Let $a\geq 2$. Taking into account the considered Case 2, we may assume that $d^-(x_m,\{x_1,x_2, \ldots , x_{a-1}\})=0$.
It is not difficult to see that  $$A(\{x_1,x_2,\ldots , x_{a-1}\}\rightarrow 
\{x_{a+1},x_{a+2}\})=\emptyset.$$
Indeed,  if $x_ix_j\in A(D)$ with $i\in [1,a-1]$ and $j\in [a+1,a+2]$, then $C(z)=x_1x_2\ldots x_ix_jx_{a+2}\\Rx_{i+1}\ldots x_ax_mzx_1$ is a cycle of length at least $m+2$, a contradiction.
Thus, we may assume that 
$$
 A(\{x_1,x_2,\ldots , x_{a-1}\}\rightarrow 
\{x_{a+1},x_{a+2},x_{a+3}= x_{m}\})=\emptyset. \eqno (14)
 $$ 
Since $D-x_{a}$ is strong, it follows that in    $D-x_{a}$ there is a path from a vertex  $x_i\in \{x_1,x_2,\ldots , x_{a-1}\}$ to a vertex 
$x_j\in \{x_{a+1},x_{a+2}, x_{a+3}=x_m\}$. Let $Q$ be such a  shortest  path. Then, using (2) and $d(z, Y\cup \{x_{a+1},x_{a+2}\})=0$, it is not difficult to see that $x_j\in \{x_{a+1},x_{a+2}\}$ and the internal vertices of $Q$ are in $\sum_{i=d+1}^{f}V(H_i)$.  
This means that $C_{m+1}(z)$ contains an $(x_i,x_j)$-bypass such that $z\notin V(C_{m+1}(z)[x_{i+1},x_{j-1}])$. Therefore by Lemma 3, $D$ contains a $C(z)$-cycle of length at least $m+2$, a contradiction.
Let now $a=1$. Then $b=l+1=m=4$.  From 
$d^-(z)\geq 2$, $d^-(z)\geq 2$ and $d(z, Y\cup \{x_2,x_3\})=0$ it follows that 
 $x_1\rightarrow z\rightarrow  x_4$, which contradicts Claim 1(i). The discassion of Case 3 is completed.

 \textbf{Case 4.} $a=l-1$. 

Taking into account  Cases 2-3 and the digraph duality, we may assume that $b=l+1$ and $A(\{x_{1},x_{2},\ldots , x_{l-1}\}\rightarrow  \{x_{l+2},x_{l+3},\ldots , x_{m}\})=\emptyset$.

Assume first that $m\geq b+1=l+2$. If $x_l\rightarrow x_i$ with $i\in [l+2,m]$, then the cycle $C(z)=x_1x_2\ldots x_{l-1}x_{l+1}\ldots x_{i-1}Rx_lx_i\ldots x_mzx_1$ has length at least $m+3$, a contradiction. We may therefore assume that $d^+(x_l,\{x_{l+2},\ldots , x_m\})=0$. As a result, we have
$A(\{x_{1},x_{2},\ldots , x_{l}\}\rightarrow  \{x_{l+2},x_{l+3},\ldots , x_{m}\})=\emptyset$, which contradicts Proposition 1.

Assume second that $m=b=l+1$.
Let $a\geq 2$. Taking into account the considered cases,
 it is not dificult to show that $A(\{x_{1},x_{2},\ldots , x_{a-1}\}\rightarrow  \{x_{a+1},x_{a+2}=x_{m}\})=0$. Since $D-x_a$ is strong, there is a path from a vertex $x_i\in \{x_1,x_2,\ldots , x_{a-1}\}$  to a vertex $x_j\in \{x_{a+1}, x_{a+2}\}$. Let $Q$ be such a shortes path. From Claime 1(i), (2) and $d(z,\{x_{a+1}\}\cup Y)=0$ it follows that $Q$ is a $(C(z),x_i,x_l)$-bypass, whos internal vertices are in $\cup_{i=d+1}^fV(H_i)$ and
$z\notin V(C_{m+1}(z)[x_{i+1},x_{l-1}])$, this contradicts Lemma 3.
Let now $a=1$. Then $m=3$ and $x_1\rightarrow z\rightarrow x_3$, which contradicts Claim 1(i). This contradiction  completes the discussion of Case 4. Lemma 4 is proved. \fbox \\\\

 Now we are ready to prove the main result of this note. 
 
\textbf{Proof of Theorem 6:} Let $D$ be a 2-strong digraph of order $n$ satisfying the conditions of THeorem 6. Suppose that $D$ is not Hamiltonian. Let $C_m(z):=x_1x_2\ldots x_mx_1$ be a longest cycle through $z$ in $D$ and let $Y=V(D)\setminus V(C_m(z))$. Note that $m\geq 3$ since $D$ is 2-strong. Recall that  by the supposition 
 of the theorem, $m\geq n-k-2$.  Since $C_m(z)$ is a longest cycle in $D$, it follows that any vertex $u\in Y$ cannot be inserted into  $C_m(z)$. Then by Lemma 1, $d(u,V(C_m(z))\leq m$ and
$$
n+k\leq d(u)=d(u,V(C_m(z)))+d(u,Y)$$ $$\leq m+2n-2m-2=2n-m-2.
$$
Hence, $m\leq n-k-2$ (i.e., $m=n-k-2$), $|Y|=k+2\geq 2$, $d(u,Y)=2k+2$ (i.e., $D\langle Y\rangle$ is a complete digraph) and $d(u,V(C_m(z)))=n-k-2$. If some vertex of $Y$ is adjacent to every vertex of $C_m(z)$, then $D$ is Hamiltonian. 
We may therefore assume that there are vertices $y\in Y$ and $x_i$, say $x_{n-k-2}$, 
which  are not adjacent. Using the facts that $d(y,V(C_m(z)))=n-k-2$, $D\langle Y\rangle$ is complete and  Lemma 2,  
it is not difficult to show that  $d(x_{n-k-2},Y)=0$, $x_{n-k-3}\rightarrow Y \rightarrow x_1$ and $z=x_{n-k-2}$.
Since for every vertex $u\in Y$, $d(u,V(C_{n-k-2}(z)))=n-k-2$,  it follows that $D$ satisfies the conditions of Lemma 4. Therefore, $d(z)=d(z,\{x_1,x_2,\ldots , x_{n-k-3}\})\leq n-k-5$, which contradicts that $d(z)\geq n-k-4$. The theorem is proved. \fbox \\\\

\section {Conclusion}

For Hamiltonicity of a graph  $G$ (undirected graph), there are
numerous sufficient conditions in terms of the connectivity  number $k(G)$ of $G$, where $k(G)\geq 3$ (recall that for a graph $G$ to be Hamiltonian, $k(G)\geq 2$ is a necessary condition) and the minimum degree $\delta(G)$ or the sum of degrees of some vertices with certain properties. Results on Hamiltonian graphs can be found in the survey papers by Gould, e.g. \cite{[15]}. 
This is not the case for the general digraphs. Moreover, in \cite{[16]}, the author proved that: For every pair of integers $k\geq 2$ and $n\geq 4k+1$ (respectively, $n=4k+1$), there exists a $k$-strong  $(n-1)$-regular (respectively,  with minimum degree at least $n-1$ and with minimum semi-degrees at least $2k-1=(n-3)/2$) a non-Hamiltonian digraph of order $n$. 

There are a number of degree or degree sum  condition for a bipartite digraph to be Hamiltonian. The reader can find more information on the topic in survey paper \cite{[17]} by  Ge,  Ye and Zhang. Often, the lower bounds in such conditions are best possible. However, many reseachers reduce the bounds and try to  identify all exceotional bipartite digraphs, that is the non-Hamiltonian digraphs satisfing these new conditions, see \cite{[18]} and the papers cited there.

Based on these and  the evidence from Theorem 4, we propose the following problem.

\textbf{Problem 2:}  {\it Investigate the Hmiltonicity of bipartite digraphs by requring  that the degree condition satisfies only for some vertices or some pairs of vertices with an additional restriction (for detailes, see arXiv:2306.16826)}.

\end{document}